# Evaluating Pricing Strategy Using e-Commerce Data: Evidence and Estimation Challenges

**Anindya Ghose and Arun Sundararajan**

*Abstract.* As Internet-based commerce becomes increasingly widespread, large data sets about the demand for and pricing of a wide variety of products become available. These present exciting new opportunities for empirical economic and business research, but also raise new statistical issues and challenges. In this article, we summarize research that aims to assess the optimality of price discrimination in the software industry using a large e-commerce panel data set gathered from Amazon.com. We describe the key parameters that relate to demand and cost that must be reliably estimated to accomplish this research successfully, and we outline our approach to estimating these parameters. This includes a method for "reverse engineering" actual demand levels from the sales ranks reported by Amazon, and approaches to estimating demand elasticity, variable costs and the optimality of pricing choices directly from publicly available e-commerce data. Our analysis raises many new challenges to the reliable statistical analysis of e-commerce data and we conclude with a brief summary of some salient ones.

## 1. INTRODUCTION

The adoption of Internet-based commerce has provided academic researchers with a wealth of new data on demand and pricing across a number of industries. The availability of these data and their growing use in empirical studies of electronic commerce raises a number of new statistical and econometric issues. In this article, we describe how we empirically analyze and evaluate pricing strategy in the consumer software industry using a large-scale e-commerce data set from Amazon.com. We describe some of the methods we have applied to our analysis of these data, how we have adapted them to address issues unique to e-commerce data, and we summarize open challenges whose resolution will help facilitate more robust empirical research in electronic commerce.

Pricing strategy in the consumer software industry (and in many other industries) often involves the use of *price discrimination*, which, broadly, aims to identify (directly or otherwise) customers who are willing to pay more for a product and to charge them a higher price. Beyond the notion of "first-degree" price discrimination which involves charging different consumers different prices for an identical good (Aron, Sundararajan and Viswanathan, 2006;

*Anindya Ghose is Assistant Professor of Information, Operations and Management Sciences, and Arun Sundararajan is Assistant Professor of Information Systems and Director (IT Economics) of the Center for Digital Economy Research, Leonard N. Stern School of Business, 44 West 4th Street, New York, New York 10012, USA e-mail: aghose@stern.nyu.edu; asundara@stern.nyu.edu.*

Choudhary et al., 2005), there are a variety of ways that firms price-discriminate. For example, a seller may price differently depending on whether a consumer has purchased from the firm before (these are typically called *introductory offers*). A seller may vary the price of a product depending on how many units of the product are purchased by an individual consumer; this is commonly referred to as *non-linear pricing* (Sundararajan, 2004b). A seller may base the price of a product on whether other related products are also purchased from the same firm: this is called *bundling* (Bakos and Brynjolfsson, 1999), and a seller may choose to implement either *pure* bundling, under which a set of products are sold only as a bundle, or *mixed* bundling, under which both the bundle and individual products are sold (Ghose and Sundararajan, 2005b). As an example of the latter, Microsoft sells its Office suite of software as a bundle of Word, Excel and PowerPoint in addition to selling each of these products individually. A seller may create different but related versions of a product (typically one of higher quality or with more features) and price them differently. This is referred to as *versioning* and aims to price-discriminate by exploiting differences in how much different customers value product quality. There are multiple versions of a large number of popular desktop software titles that differ only in their quality or number of features (rather than in their development or release date) and that are sold at different prices. Current examples include Adobe Acrobat, TurboTax, Microsoft Money and Norton AntiVirus. These are examples of software titles for which a firm has developed a flagship version, disabled a subset of the features or modules of this version, and released both the higher quality version and one or more lower quality versions simultaneously. A related form of price discrimination is based on releasing *successive generations* of the same product in multiple periods, with a period of time where the old and new generations overlap; since each new generation represents an improvement in the overall performance of the product, the simultaneous presence of two or more successive generations is analogous to the presence of two or more related products of varying quality (Ghose, Huang and Sundararajan, 2005).

The objective of a software company that price-discriminates is to maximize the profits it generates from the sale of its products. However, price discrimination can often have countervailing effects on a firm's profits. For instance, two consequences of introducing a lower quality version of an existing product to price-discriminate are the loss of profits from customers who switch from purchasing the higher quality version to purchasing the lower quality version (commonly termed *cannibalization*) and a gain in profits from new customers, for the lower quality version, who either did not purchase the product earlier or who purchased a competing product. [In many ways this is similar to the cannibalization that occurs when used products compete simultaneously with new products (Ghose, Telang and Krishnan, 2005; Ghose, Smith and Telang, 2006).] The interplay between these consequences eventually determines the optimality of versioning. A similar pair of consequences, with opposing effects, characterizes the eventual profitability of bundling. Similarly, nonlinear pricing that discounts high usage levels too extensively can reduce a seller's profits.

Thus, to profit from price discrimination, a software company must make an appropriate choice of the form of price discrimination; it must choose its prices optimally and sometimes it must determine optimal quality levels for an inferior (related) set of products or the size of a bundle. There is no published research with evidence that software companies in fact make these price-discrimination choices optimally; however, the availability of detailed price and demand data from e-commerce sites like Amazon now makes it feasible to empirically assess the optimality of their choices. A first goal of our research program is therefore to use these data to evaluate the optimality of such price-discrimination strategies in the software industry empirically. This is a problem of significant economic importance.

To do so, one first needs a method to convert "sales ranks" reported by Amazon.com into actual demand levels. Amazon publishes a sales rank for each product it sells, which is the rank of the product within its category based on recent demand (more on this later). Next, the demand system associated with our products (i.e., how the variation in prices is associated with variation in demand) needs to be estimated. Amazon.com does not provide any data about the variable cost of the products it sells: we therefore also need to infer these costs from our data (since the profit to a seller is determined not just by price charged and quantity sold, but also by its cost per unit). We describe our approach to accomplish each of these goals. We briefly summarize other estimates that contribute to our research program and we conclude with some of the key statistical challenges that emerge from our analysis. The flow chart in Figure 1 indicates the steps we prescribe for estimating the optimality of pricing.

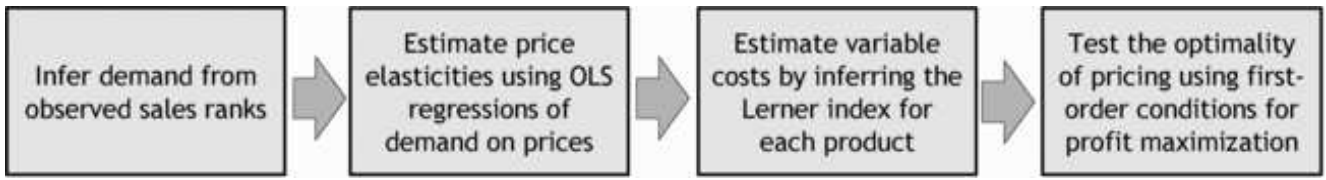

Fig. 1. *Sequence of steps to determine the optimality of software pricing from e-commerce data.*

## 2. SUMMARY OF DATA

Our data are compiled from publicly available information on new software prices and sales rankings at Amazon.com, the largest on-line retailer of consumer software. Our data are gathered using automated Java programs to access and parse HTML and XML pages downloaded from its web site, three times each day, at equally spaced intervals. Our sample contains 330 products randomly selected from each of four major categories—business and productivity, security and utilities, graphics and development and operating systems software. (Our random sample was created by first compiling a list of software products that were sold on Amazon during the year and then using Excel's random number generator to choose from them. We chose a sample size of 330 since this yielded what we felt was a sufficient number of distinct titles within each major category.) We collect all relevant data on list prices (the manufacturer's suggested price), new prices (the price charged by Amazon.com), sales ranks (to be discussed further soon), product release date, average customer review and number of reviewers who contributed to this average. To facilitate a clearer understanding of how each of these pieces of information is reported to a consumer on Amazon.com's web site, a screen shot of an Amazon page is illustrated in Figure 2.

Fig. 2. *How the data we gather from Amazon.com are displayed on its web site.*

Table 1
*Summary statistics*

| Variable | Mean | Std. dev. | Min | Max |
|---|---|---|---|---|
| *Sales rank* | 1649.61 | 1971.26 | 1 | 11622 |
| *List price* | 69.16 | 226.17 | 19.95 | 1799.99 |
| *Amazon price* | 65.53 | 208.57 | 14.95 | 1699.99 |
| *New non-Amazon price* | 17.74 | 23.08 | 10.01 | 209.99 |
| *Customer rating* | 3.14 | 0.99 | 1 | 5 |
| *Number of reviewers* | 25.72 | 66.3 | 1 | 606 |
| *Days release* | 717.7 | 1336.22 | 0 | 1750 |

Our sample consists of products that have different versions as well as products that are sold as bundles (in addition to the individual components). We are able to determine this because software manufacturers use terms like "premier," "deluxe" or "standard" to denote versions of the same title that vary in quality (which is typically measured by the number of features). Similarly, a product suite that contains individual products has the term "bundle" associated with it. The details of individual components within each bundle are provided in a bundle's product description. Based on the release date of the product, we can infer if it is the current or the previous year's edition.

We also collect data on secondary market activity, including used prices (prices charged by sellers who have posted second-hand copies of the product for sale) and new prices from non-Amazon sellers (these are sellers who are not affiliated with Amazon but are allowed to sell goods on Amazon in exchange for a commission on the transaction price). Table 1 provides summary statistics of our data.

We have categorized our software titles in three ways: (i) based on those titles that have just two versions and those that have more than two versions, (ii) based on whether the title is sold as part of a bundle of other products or as a stand-alone product, and (iii) based on whether the title is from the most recent generation or from a previous generation. This categorization is summarized in Table 2. Thus, for example, our sample contains 32 unique titles which each have two versions (a higher quality and a lower quality version). Similarly, our sample contains 56 unique titles which each have both the current and the previous generation available simultaneously. The other rows can be interpreted in a similar way.

Our data were collected between January 2005 and November 2005. (For the duration of our study there were a few instances during which the Java program was unable to collect data all three times during the day. In most cases this happened if the Amazon server was not functioning properly during the time the data were being gathered. However, this does not affect our analysis, primarily because of the low frequency with which prices are changed by Amazon. In general, we find that there are far fewer price changes than changes in sales rank; thus any missing information on prices within the same day will have almost no impact on our estimates of price elasticities.) The distributions of sales ranks and retail (Amazon) prices across our products are summarized in Figure 3(a) and (b). We also provide a scatterplot of prices and sales ranks in Figure 3(c).

## 3. ESTIMATION AND PRICING

### 3.1 Demand Estimation

Amazon.com does not report its periodic demand levels. Instead, it reports a sales rank for each product sold on its site, which ranks the demand for a product relative to other products in its category. Thus, the lower the cardinal value of the sales rank, the higher the demand for that particular item. Prior research (e.g., Chevalier and Goolsbee, 2003; Brynjolfsson, Hu and Smith, 2003) has associated these sales ranks with demand levels. To do so, the authors assume that the rank data have a Pareto distribution (i.e., a power law). They then convert sales ranks into periodic demand levels by conjecturing the Pareto relationship $\log[Q] = \alpha + \beta \log[rank]$, where $Q$ is the (unobserved) demand for a product, $rank$ is the (observed) sales rank of the product and $\alpha, \beta$ are industry or category-specific parameters. [Chevalier and Goolsbee (2003) reported that evidence that the Pareto distribution fits well can be found using the weekly *Wall Street Journal* book sales index, which, unlike other bestseller lists, gives an index of the actual quantity sold. This index is constructed by surveying Amazon.com, BN.com,

Table 2
*Various product categories in sample*

| Product category | Number of unique titles | Total number of products |
|---|---|---|
| Bundles | 68 | 136 |
| Versions (2) | 32 | 64 |
| Versions (>2) | 19 | 57 |
| Successive generation | 56 | 112 |

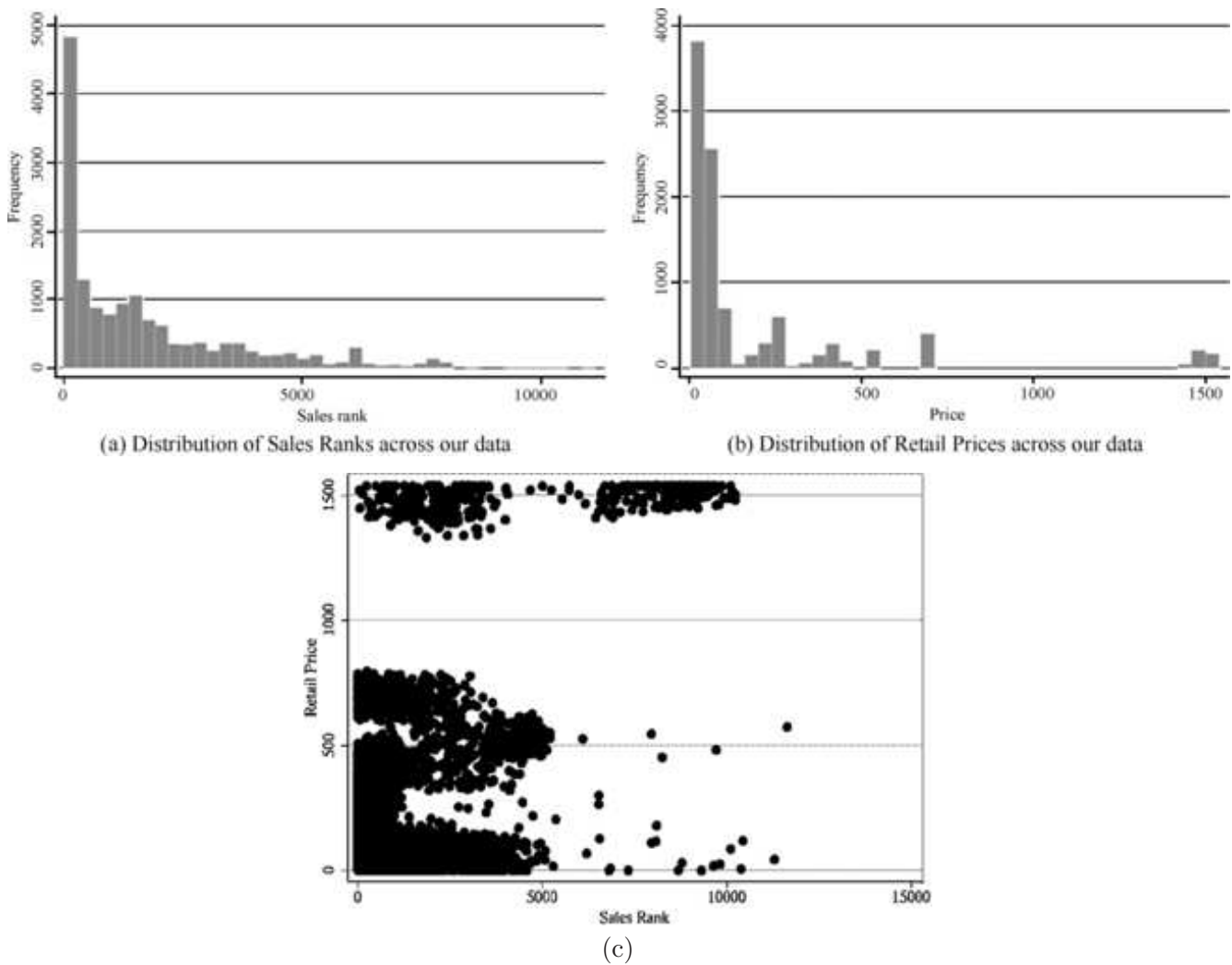

(a) Distribution of Sales Ranks across our data

(b) Distribution of Retail Prices across our data

(c)

Fig. 3. *The distribution of* (a) *sales ranks and* (b) *prices across observations in our data set. The histograms are based on a single entry per product.* (c) *Scatterplot of retail prices and sales ranks for all products in our sample.*

and several large brick and mortar book chains. For discussions on the use of power law distributions to describe rank data, see Pareto (1896/1897) and Quandt (1964).]

A number of recent studies (e.g., Ghose, Smith and Telang, 2006) pertaining to the book industry have used estimates of $\alpha$ and $\beta$ from this prior literature. However, these are industry-specific parameters and, to our knowledge, there are no corresponding estimates available for software. Furthermore, in summer 2004, Amazon altered its sales rank system in the following way: it eliminated its three-tier system, updating ranks each hour for most products (rather than merely for the top products), and it moved to a system that uses exponential decays to give more weight in the sales rank to newer purchases. The exact details of the calculation are pro-prietary to Amazon (e.g., the half-life of the decay). In the earlier three-tier system, there were three distinct ranking schemes on Amazon: one for the top selling 10,000 products, another for the products between 10,000 and 100,000, and a third for ranks above 100,000. Products with sales ranks between 1 and 10,000 were reranked every hour, products in the range from 10,000 to 100,000 were reranked once a day and products with ranks greater than a 100,000 were updated once a month (Chevalier and Goolsbee, 2003). The current system involves reranking all products every hour.

Toward a more current and accurate reverse engineering of the ranking system to infer actual periodic demand, we have designed and conducted an independent analysis to convert measured sales ranks into demand data. Retaining the assumption of a Pareto relationship between demand and sales rank, we combine a series of purchase experiments with the analysis of a time series of sales ranks of all the 330 products in our sample to estimate both $\alpha$ and $\beta$.

Our purchasing experiment proceeded as follows. Over a two-week period in mid-June 2005, we collected hourly sales rank data for each of the products in our panel, yielding a time series of 336 observations for each product. For products not ranked too high, a general trend in these time series is an extended downward drift in the rank value over many hours (i.e., the rank becomes a progressively larger number), followed by intermittent spikes which result in a large upward shift in rank (i.e., the rank became a smaller number suddenly). This is illustrated for two candidate products in Figure 4. We interpret these spikes as reflecting time periods in which one or more purchases have occurred.

This procedure yielded a data set of a certain number of observations, which associated a weekly demand level with each average sales rank, for two successive weeks. Weekly unit sales ranged from 0 to 16. Using the implied pairs of average weekly demand and average sales rank, we then estimated the ordinary least squares (OLS) equation

$$(1) \qquad \log[Q+1] = \log[\alpha] + \beta \log[sales\ rank],$$

where $Q$ is average weekly demand, and *sales rank* is the corresponding average sales rank. [Similar to Brynjolfsson, Hu and Smith (2003), we used White's heteroskedasticity-consistent estimator (see Greene, 2000, page 463) to estimate both parameters.] The results of these experiments yielded $\alpha = 8.352$ and $\beta = -0.828$. The following list provides a sense of what these estimates imply:

Weekly sales of two units correspond to an average sales rank of about 3100.

Weekly sales of 10 units correspond to an average sales rank of about 440.

Weekly sales of 25 units correspond to an average sales rank of about 150. [The standard errors for $\alpha$ and $\beta$ were 0.042 and 0.032, resp., and the estimates were significant at the 1% level. Further details of this experiment and its results are presented in Ghose and Sundararajan (2005a).]

An interesting aspect of our approach is that it allows one to characterize a number of economic measures of interest (price markups and demand elasticities) based purely on sales rank ratios and the parameter $\beta$. This provides a framework for a wider range of future empirical research in e-commerce and also reduces the extent to which one's results are affected by the error in estimating $\alpha$.

Our sample encompasses multiple products with observations collected over time, and our data set therefore has elements of both cross-sectional and time-series data. Consistent with existing published and current research, we analyze our observations as panel data (for a detailed treatment of the econometric analysis of panel data, see Wooldridge, 2002).

### 3.2 Estimates of Price Elasticity

Given the price variation across products and across time, and the measures of quantity in each period that come from the sales ranks, we can infer the price sensitivity of consumers to on-line product sales. This would require estimating own- and cross-price elasticities of the products in our sample. The price elasticity of demand is a measure of the sensitivity of demand to price changes. Specifically, the own-price elasticity of demand is calculated as the percentage change in demand caused by a unit percentage change in a product's own price, and the cross-price elasticity of demand is calculated as the percentage change in demand caused by a unit percentage change in another product's price. In our context, the other product could be either a different version (high or low quality) of the same good or a component of a bundle. Own-price elasticities are generally negative—the quantity of a product sold decreases as its price increases. On the other hand, cross-price elasticities can be either positive or negative. If $X$ and $Y$ are substitute goods (e.g., the two versions of a product), the cross-price elasticity of demand is positive; that is, the quantity of good $X$ varies directly with a change in the price of good $Y$. If $X$ and $Y$ are complementary goods (e.g., computer hardware and software), the cross-price elasticity of demand may be negative; that is, the quantity of good $X$ varies inversely with a change in the price of good $Y$. Figure 5 shows a plot of the changes in prices and sales ranks for the two versions (of high and low quality) of a specific product for a specific period of time.

To compute own-price and cross-price elasticities, we estimate OLS regressions which control for unobserved heterogeneity across products and across categories (we use the fixed effects transformation). Based on these estimates, we subsequently compute the own- and cross-price elasticities by weighting them with the appropriate Pareto mapping parameter $\beta$, which was estimated earlier. In other words, estimating the equations using log ranks, rather than actual quantities, yields the correct elasticities, but they are scaled up by the Pareto mapping parameter. This is similar to the approach used in prior literature (Chevalier and Goolsbee, 2003; Ghose, Smith and Telang, 2006). These regressions have the general form

$$\log(rank_{it}) = a + \phi \log(p_{it}) + \sum_{j \in S_i} \gamma_j \log(p_{jt})$$

$$(2) \qquad\qquad + \lambda \log(\hat{p}_{it}) + \omega X + \varepsilon_{it},$$

where $i$ indexes the product in question (e.g., the high-quality version of a specific title), $t$ indexes the date, $S_i$ is the set of products whose prices affect the demand for product $i$ (e.g., the price of a lower quality version corresponding to a high-quality version or the price of a bundle which contains product $i$ as one component), $\hat{p}_{it}$ is the lowest price posted for the corresponding non-Amazon marketplace product (the best price across all conditions by competing sellers on Amazon's secondary market) and $X$ is a vector of control variables. Our control variables include the time since the product was released (*days release*), the average customer rating (*customer rating*) and the number of reviewers (*number of reviewers*) who have reviewed the product. [We use the fixed effects (within) transformation (Wooldridge, 2002, Chapter 10.5) to control for unobserved heterogeneity across products.]

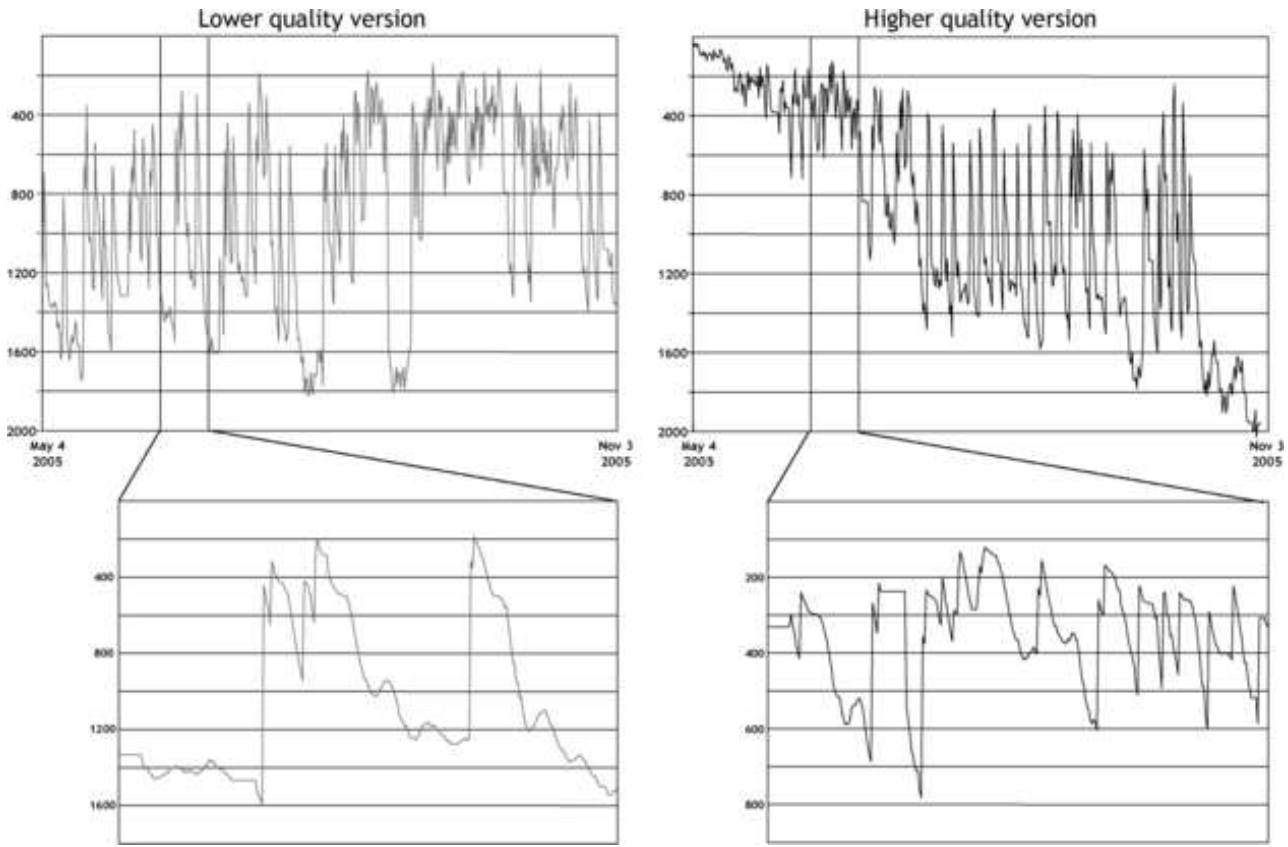

Fig. 4. *The variation of sales rank with time for the higher and lower quality versions of one of the software titles we track. The charts on top illustrate sales ranks gathered during successive 8-hour periods, charted for a 6-month window of our data. The charts below graph hourly sales rank data for the specific shorter time window (between the vertical lines on the upper charts) and illustrate the sales rank spikes associated with sales. The flat portions of these lower charts reflect (short) intervals where Amazon did not update the sales ranks that it published on its web site.*

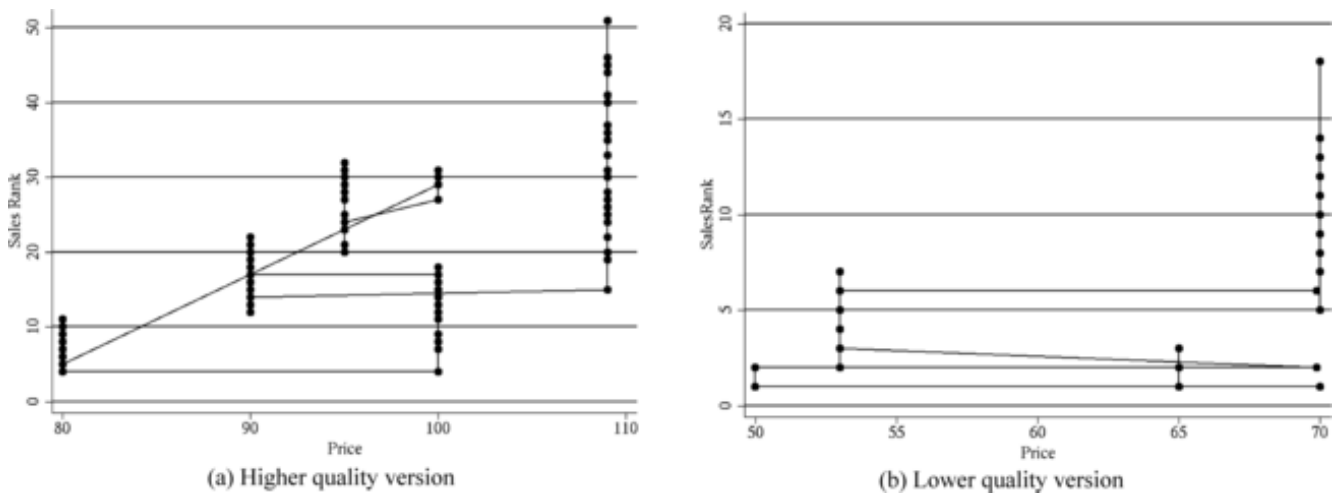

Fig. 5. *The variation of sales rank with retail price for two versions of a specific software title in our data set. A line between two points indicates that the ( price, sales rank) had changed from one of the points to the other in successive periods.*

We can use the results of this regression to calculate the relevant own- and cross-price elasticities $\eta_{ii}$ and $\eta_{ij}$, respectively. Note that $\phi$ is a measure of how sensitive the product's own sales rank is to the product's own price, while $\gamma_j$ is a measure of how sensitive the product's own sales rank is to the price of the substitute product. Let $(p_i, Q_i)$ represent the quantity and price of product $i$ and let $(p_j, Q_j)$ represent the quantity and price of product $j$. One can easily show that having estimated the parameters $\phi$ and each of the $\gamma_j$'s from the above regression, the own-price elasticity of demand for product $i$ (Mas-Collel, Whinston and Green, 1995) is given by

$$(3) \qquad\qquad \eta_{ii} = \beta\phi = \frac{\partial Q_i}{\partial p_i} \times \frac{p_i}{Q_i},$$

where $\beta$ was estimated in Section 3.1 and the cross-price elasticity of demand for product $i$ with respect to product $j$ is

$$(4) \qquad\qquad \eta_{ij} = \beta\gamma_j = \frac{\partial Q_i}{\partial p_j} \times \frac{p_j}{Q_i}.$$

These elasticity estimates describe how demand varies with price, and form the basis for analyzing the optimality of a firm's chosen price-discrimination strategy, since they enable us, for example, to assess how demand would vary if the firm altered its price discrimination by removing a version or discontinuing a bundle. They also are inputs to the estimation process for variable costs, as described in the following section. (In contrast to some other demand estimation models, we do not have data on other inputs to the marketing mix, such as advertising. It is conceivable that sales ranks may also be affected by off-line advertising. In the absence of data, we are unable to capture this effect in our model.)

### 3.3 Cost Estimation

Many products in information technology (IT) industries have an unusual cost structure: high fixed costs of production, but near-zero or zero variable costs of production. This cost structure characterizes a class of technology products which are collectively termed *information goods*. Put differently, the cost of producing the first unit of an information good is very high, yet the cost of producing each additional unit is virtually nothing. For instance, Microsoft spends hundreds of millions of dollars to develop each version of its Windows operating system. Once this first copy of the operating system has been developed, however, it can be replicated costlessly, which leads to widespread piracy, a factor which can be incorporated into pricing following Sundararajan (2004a), but which we do not explicitly model in this study. Early examples of information goods were computer-based information services and software; currently, a wide variety of diverse products—video, music, textbooks, digital art, to name a few—share this unique cost structure.

Contrary to what is commonly assumed in the IT economics literature, packaged consumer software is not an "information good." It has positive variable costs associated with its production, packaging and distribution, and these may represent a substantial fraction of the price of such software, especially since a number of titles are priced under fifty dollars. Therefore, to assess the optimality of a seller's choice of price discrimination, we need estimates of the variable costs of the software titles in our data set. We estimate the variable costs by inferring the Lerner index for each product version $i$, defined as the ratio of the markup to the price, that is, $((p_i - c_i)/p_i)$, where $p_i$ is the retail price and $c_i$ is the variable cost of product $i$. Markup is simply defined as price minus marginal cost, that is, the margin on each unit of sale.

To do this estimation reliably using e-commerce data, we have developed an extension of the model of Hausman (1994) that provides a way to estimate markups using just sales rank data and prices. We begin with the approach of Hausman (1994), who provides the following equation to estimate the markups for products sold by multiproduct oligopolists, weighted by their market share. Since software firms generally sell multiple products and compete with multiple firms in the market, it is important to consider this formulation of the Lerner index. Consider a set of related products indexed by $i$. The first-order conditions for oligopoly profit maximization yield the system of equations

$$
s_j + \sum_i \left[ \left( \frac{p_i - c_i}{p_i} \right) s_i \right] \eta_{ij} = 0,
\tag{5}
$$

$$
j = 1, 2, \ldots, n.
$$

Here, $s_i$ is the demand share of product $i$ (demand share is the ratio of revenues from product $i$ to the total revenues from all related products), $\eta_{ii}$ is product $i$'s elasticity of demand with respect to its own price and $\eta_{ij}$ is the cross-price elasticity of demand with respect to the price of product $j$. We therefore have a system of linear equations

$$
s + N'm = 0,
\tag{6}
$$

where $s$ is the vector of revenue shares, $N$ is the matrix of cross-price elasticities $[\eta_{ij}]$ and $m = [m_0, m_1, \ldots, m_n]$, where

$$
m_i = \left( \frac{p_i - c_i}{p_i} \right) s_i
\tag{7}
$$

is the Lerner index of product $i$ multiplied by its product share. The marginal costs $c_i$ of each individual product can then be estimated by inverting $N$ to solve the system of equations (6).

Our extension of this approach allows the estimation of variable costs using just sales ranks, the parameter $\beta$ which we estimate in Section 3.1, and observed retail prices. Our system of equations for a set of related products $0, 1, \ldots, n$ with prices $p_i$ and sales ranks $R_i$ is derived from the above equations as

$$
s_j = -\sum_i \left[ \frac{p_i - c_i}{p_i} s_i \right] \eta_{ij}, \quad j = 1, 2, \ldots, n,
\tag{8}
$$

which implies that

$$
s_j = \sum_i \left[ \frac{p_i - c_i}{p_i} s_i \right] \beta \gamma_j,
\tag{9}
$$

which in turn implies that

$$
s_j = \sum_i \left[ \frac{p_i - c_i}{p_i} \right] \beta \left( \frac{p_i}{R_j} \frac{dR_j}{dp_i} \right) s_i
\tag{10}
$$

or

$$
s_j = \beta \sum_i \left[ (p_i - c_i) \frac{s_i}{R_j} \frac{dR_j}{dp_i} \right],
\tag{11}
$$

where

$$
\frac{1}{s_i} = 1 + \frac{p_i}{p_j} \left( \frac{R_j}{R_i} \right)^\beta.
\tag{12}
$$

### 3.4 Optimality of Pricing

In this section, we summarize how we can test the optimality of pricing strategies by software manufacturers. Consider the case when the software firm is producing two versions of the product—a high-quality version and a low-quality version. The total profit from a pair of versions $i$ and $j$ is thus

$$
\pi = k(p_i - c_i)Q_i + k(p_j - c_j)Q_j.
\tag{13}
$$

First-order conditions for profit maximization with respect to prices yield the partial derivatives

$$
\frac{\partial \pi}{\partial p_i} = kQ_i + k(p_i - c_i)\frac{\partial Q_i}{\partial p_i}
$$

$$
+ k(p_j - c_j)\frac{\partial Q_j}{\partial p_i},
\tag{14}
$$

$$
\frac{\partial \pi}{\partial p_j} = k(p_i - c_i)\frac{\partial Q_i}{\partial p_j}
$$

$$
+ kQ_j + k(p_j - c_j)\frac{\partial Q_j}{\partial p_j}.
\tag{15}
$$

If the products are priced optimally, these partial derivatives should be equal to zero. Note that we can evaluate each term on the right-hand side in the above equations empirically based on our data set and the intermediate steps described in Sections 3.1, 3.2 and 3.3. Specifically, from the data we can infer what the prices $(p_i)$ as well as what the quantities $(Q_i)$ are. Based on the cost estimation procedure outlined above, we can infer the marginal costs $(c_i)$ of each product. Finally, from the demand estimation procedure, we can derive the price elasticities and, consequently, impute what the specific derivatives $\partial Q_i / \partial p_i$ are. Thus, based on the signs of these partial derivatives, we can empirically test if the firm's prices are optimal, underpriced or overpriced.

| Variable | Estimates (standard error) |
| --- | --- |
| Constant | $-0.05(1.58)$ |
| $\ln(p_{bundle})$ | $2.22^{***}(0.48)$ |
| $\ln(p_{component1})$ | $-0.19^{***}(0.06)$ |
| $\ln(p_{component2})$ | $-0.12^{***}(0.03)$ |
| $\ln(\hat{p}_{bundle})$ | $-0.24^{***}(0.07)$ |
| $\ln(days\ release)$ | $0.18^{*}(0.1)$ |
| $R^2$ | $0.42$ |

[a]Standard errors are given in parentheses. The dependent variable is $\ln(sales\ rank)$ of the bundle.

***, ** and * denote significance at the 0.01, 0.05 and 0.1 levels, respectively.

### 3.5 Examples

In this section we provide the parameter estimates for two products in our sample: Adobe Photoshop and Premier bundle, and Microsoft Office (our results are quite robust: the estimates from bootstrapping with different repetitions are the same in magnitude and direction as the original estimates), summarized in Tables 3 and 4.

Note that, as expected, the sign of the own-price elasticity is positive while the signs of the cross-price elasticities are negative (recall that an increase in sales rank implies a decrease in sales). The other control variable (days release) suggests that, as expected, sales of products decrease over time. The numbers indicate that own-price elasticity is significantly higher than cross-price elasticities for the Adobe bundle with respect to each of its two components ($p_{component1}$ and $p_{component2}$). Interestingly, we find that the cross-price elasticity of the high-

| Variable | Estimates (standard error) |
| --- | --- |
| Constant | $7.77(7.21)$ |
| $\ln(p_{professional})$ | $1.91^{***}(0.58)$ |
| $\ln(p_{standard})$ | $-2.54^{***}(-0.97)$ |
| $\ln(\hat{p}_{professional})$ | $-0.36^{***}(-0.11)$ |
| $\ln(days\ release)$ | $0.01^{***}(0.003)$ |
| $R^2$ | $0.32$ |

[a]Standard errors are given in parentheses. The dependent variable is $\ln(sales\ rank)$ of the high-quality version which is Microsoft Office Professional.

***, ** and * denote significance at the 0.01, 0.05 and 0.1 levels, respectively.

quality version of Microsoft Office with respect to the low-quality version ($p_{standard}$), is actually higher than the own-price elasticity. This highlights that consumer demand for Microsoft Office Professional is very sensitive to the price of Microsoft Office Standard. Also we do not see much variance in the extent to which competing prices from non-Amazon sellers matter in influencing demand at Amazon; in both cases, the cross-price elasticities (from competing sellers) are significantly lower than the own-price elasticities. [These parameter estimates are obtained from OLS regression models of the kind mentioned earlier. Also, because of the structure of this industry, quantity and price are not jointly determined; thus we do not face the endogeneity concerns that would normally arise in demand regressions. With regard to Amazon's own price, because software titles are produced in large quantities prior to going to market, the quantity of new products Amazon can sell is predetermined (and usually virtually infinite) at the time Amazon sets its price. This follows similar approaches taken in the literature for demand estimation of Internet product sales using the above assumed functional form for demand (Ghose, Smith and Telang, 2006; Ghose and Sundararajan, 2005b; Ghose, Huang and Sundararajan, 2005).]

As an example of a test for the optimality of pricing, we take the case of Microsoft Office. Using the estimates for own- and cross-price elasticities derived for both the high- and low-quality versions (see, e.g., Table 4 which reports the estimates for the high-quality version) in equations (14) and (15), we find that the estimated derivative of profits with respect to $p_{professional}$ is $-75.2$ and with respect to $p_{standard}$ is $-334.9$. The actual magnitudes of these estimates do not lend themselves easily to interpretation. However, their signs suggest that both versions of Office are overpriced, since they are priced at a point where the slope of the profit function is negative. [These estimates are based on maximizing the profits of the channel as a whole (i.e., the sum of the profits of the retailer and the software manufacturer). Since we do not separate the optimization problems of each of these firms, our estimates do not identify whose actions need to be changed to rectify this mispricing.]

## 4. CONCLUSION

Our objective in this paper is to outline analyzing the optimality of pricing strategy of software firms using e-commerce panel data. While we have shown how the widespread availability of e-commerce data presents a number of novel empirical research opportunities, it is important to point out that there are significant new challenges faced by researchers who aim to analyze these data in a statistically valid and economically meaningful way. Although our context is price discrimination in software, the methods we use apply equally well to e-commerce data about any consumer product category.

A key statistical challenge in demand estimation of this kind is that the time structure of e-commerce data is not well understood. Granted, one can control for systematic seasonal effects (such as time of the day or month of the year that the data were collected) and for major event effects (such as the release of a new version of Windows), and check one's data for autocorrelation. However, e-commerce is still at a relatively early stage of its evolution and the fraction of retail demand fulfilled by e-commerce sites continues to grow over time. This is driven by an increase in both the number of consumers who shop on-line and in the fraction of their purchases made on-line. Each of these factors may affect the relationship of observed e-commerce demand and price, which in turn suggests that e-commerce data may have a complex underlying time structure.

Furthermore, new theory that models the time structure of such e-commerce data in a more precise way, and techniques that identify and account for time variation, may enable future research to assess whether the demand process that generates such observations is stationary and whether the e-commerce market in question is, in fact, in equilibrium. This is a challenge not just for retail panel data, but for other forms of data generated by consumers who interact with e-commerce sites, such as bidding/reputation data from on-line auctions. Current research that studies the time structure of bid paths on eBay (Bapna, Jank and Shmueli, 2004; Jank and Shmueli, 2006) may be a first step toward understanding similar data generation processes.

A different challenge relates to the extent to which one can conclude that inferences from data sets such as ours are representative of the characteristics of an industry (in our case, consumer software) in general. Clearly, this is likely to be less of an issue as a larger fraction of commerce is conducted electronically. We have benchmarked our price and demand distributions with a comparable data set from Buy.com, another large software retailer. However, the frequency with which the latter site updates its sales ranks is different from that of Amazon, and statistical techniques that enable one to assess how representative our intraday data are based on benchmark data with a different granularity would be helpful.

In addition to the demand and cost estimates we have described in this paper, our research program also involves developing econometric estimates of how consumers perceive the relative quality levels of related products and compare them to estimates based on self-reported quality assessments from Amazon.com and subjective assessments by CNET editors. Since aggregate customer feedback measures from eBay, Amazon.com and various other review sites are frequently used in e-commerce research as measures of some form of quality, statistical techniques that facilitate assigning appropriate cardinal values to e-commerce ratings data generated by consumers and editors would contribute to the foundations of this line of research. [The details of this study are available in Ghose and Sundararajan (2005a).]

To summarize, we have described a sequence of related studies that use e-commerce panel data to evaluate the optimality of different forms of price discrimination in the software industry. By describing our data, detailing our approach to estimating some important parameters and summarizing some of the issues that researchers face when conducting such statistical analyses on e-commerce data, we have aimed to stimulate thought about statistical challenges that arise when conducting research based on these increasingly widely used data sets. We hope that this summary will encourage future work that identifies and addresses these challenges, thereby strengthening the statistical foundations of this exciting and rapidly evolving new research area.

## ACKNOWLEDGMENTS

We thank seminar participants at New York University and participants at the First Annual Symposium on Statistical Challenges in E-Commerce for their comments, Wolfgang Jank and Galit Shmueli for their detailed feedback on an earlier draft of this paper, and Rong Zheng for outstanding research assistance in data collection. This research was partially supported by a Summer 2005 grant from the NET Institute (www.NETinst.org).